%% file: main.tex
\documentclass[10pt]{amsart}
\usepackage{import}
\usepackage{./preamble/Pakete}
\usepackage{./preamble/Befehle}

\begin{document}

\import{Inhalt/}{title}
\section{Introduction}
\import{Inhalt/}{Intro.tex}
\section{Motive of the Springer Fiber}
\import{Inhalt/}{MotiveofSpr.tex}
\section{Equivariant Springer Motives}
\import{Inhalt/}{Form.tex}
%\section{Main Results}
%\import{Inhalt/}{Result.tex}
%\appendix
%\import{Inhalt/}{appendix.tex}

\bibliographystyle{amsalpha} 
\bibliography{main}

\end{document}

%% file: Inhalt/title.tex
%!TEX root = ../main.tex

% [short]{long title}
\title[Springer Motives]{Springer Motives}

\author{Jens Niklas Eberhardt}
\address{Max Planck Institute for Mathematics, Vivatsgasse 7, 
53111 Bonn, Germany
}
\email{mail@jenseberhardt.com}  % 'optional'
%\date{February, 2016}
%\subjclass[2010]{17B10, 22E46}
%\keywords{Parabolic induction, Lie algebras, category $\mathcal{O}$, Soergel modules, Koszul duality, mixed Tate motives}
\begin{abstract}
We show that  the motive of a Springer fiber is pure Tate. We then consider a category of equivariant Springer motives on the nilpotent cone and construct an equivalence to the derived category of graded modules over the graded affine Hecke algebra.
\end{abstract}
\maketitle
\setcounter{tocdepth}{1} 
\tableofcontents

%% file: Inhalt/Intro.tex
%!TEX root = ../main.tex
\subsection{Motive of the Springer Fiber}
Let $G$ be a connected reductive algebraic group over an algebraically closed field $k.$ Denote by $\nil\subset \g=\operatorname{Lie}(G)$ the associated nilpotent cone and by $\mu:\rnil\rightarrow\nil$ the Springer resolution. 
%Both $\nil$ and $\rnil$ carry an action of $G$ and the Springer resolution is equivariant with respect to this map. 
For $N\in\nil$, denote by $\B{N}=\mu\inv(N)\subset\rnil$ the Springer fiber. 

Let $\Lambda$ be some commutative ring of coefficients.
Denote by $\DM(\Spec(k), \Lambda)$ the triangulated category of Voevodsky motives over the base field $k$ with coefficients in $\Lambda$ (see \cite{MVW}).% and by $\M(\B{N})\in\DM(\Spec(k),\Lambda)$  the motive of $\B{N}.$ We will show
\begin{theorem}[Springer Fiber is Pure Tate]\label{thm:pureTate}
The motive of the Springer fiber $\M(\B{N})\in\DM(\Spec(k),\Lambda)$ is pure Tate, that is, a direct sum of Tate motives $\Lambda(n)[2n]$ for $n\geq 0$, if either $\chara k=0$ or if $p=\chara k>0$ and the following three conditions hold:
\begin{enumerate}
\item  $p$ is a good prime for every classical group appearing as a constituent of $G.$
\item  $p>3(h+1)$, where $h$ denotes the maximum of all Coxeter numbers of exceptional constituents in $G.$
\item $p$ is invertible in $\Lambda$ or $\Spec(k)$ admits resolutions of singularities.
  %\footnote{This should perhaps be replaced by: cases in which the Jacobson---Morozov theorem holds}
\end{enumerate}
\end{theorem}
\begin{remark}
(1) Springer fibers for classical groups admit an affine paving. This is shown in \cite[Theorem 3.9]{DLP} for $\chara k=0$ and generalized in \cite[Chapter 11]{JanNil} to $p=\chara k>0$ for good primes $p$. The existence of an affine paving almost immediately implies that  $\M(\B{N})$ is Tate.

(2) For exceptional groups, the existence of an affine paving is not known. However, DeConcini--Lusztig--Procesi \cite{DLP} show a slightly weaker result, namely that for $k=\C$ the Borel--Moore-homology of $\B{N}$ is torsion free, concentrated in even degrees and generated by algebraic cycles. We show how to adapt their arguments to prove that $\M(\B{N})$ is pure Tate in this case, under the assumption on $\chara(k).$

(3) The last assumption on $\Lambda$ and $p$ ensures a good behavior of motives of singular varieties in $\DM(\Spec(k),\Lambda)$, see \cite{SKA}. For example, it guarantees the existence of the localization triangle.
\end{remark}
\subsection{Equivariant Springer Motives} Now let $k=\overline{\mathbb{F}}_p.$ %where $p$ fulfills the assumptions of Theorem \ref{thm:pureTate}.
In \cite[Chapter II]{SoeWeVir} Soergel--Virk--Wendt construct a mixed version of the Bernstein--Lunts equivariant derived category using motivic sheaves. To a linear group $H$ acting on a variety $X\in \Var_k$ with finitely many orbits they associate a $\Q$-linear tensor triangulated category of \emph{$H$-equivariant orbitwise mixed Tate motives} $\MTDer{H}{X}.$
%This is a full triangulated subcategory of the category $\DerG{H}{X}$ of $H$-equivariant $\Deri$-motives on $X$ defined in \cite[Chapter I]{SoeWeVir}, where $\Deri(-)=\DA(-,\Q)$ denotes the homotopical stable algebraic derivator of \'etale motives with rational coefficients over the category of varieties $\Var_k$ over $k.$
We use this formalism to define the category of \emph{$H$-equivariant Springer motives}
$$\MTSpring{H}{\nil}\defi\langle\mu_*(\un_{\rnil})(n)\,|\, n\in \Z \rangle_{\Delta}\subset\MTDer{H}{\nil}$$
as full triangulated subcategory of $\MTDer{H}{\nil}$ generated by the \emph{Springer motive} $\mu_*(\un_{\rnil})$ and its Tate twists. Here $H$ denotes $G$ or $G\times \Gm,$ acting on $\nil$ in the natural way.

We will show
\begin{theorem}[Motivic Derived Springer Correspondence]\label{thm:motdersprin} There is an equivalence of categories
$$\MTSpring{H}{\nil}\cong \Derb(\rmodulesZ(\Chow_H^\bullet(Z,\Q),\star))$$
between the category of $H$-equivariant Springer motives and the derived category of finitely generated graded right modules over $(\Chow_H^\bullet(Z,\Q),\star).$
\end{theorem}
Here $Z=\rnil\times_{\nil}\rnil$ denotes the Steinberg variety. Furthermore, $(\Chow_H^\bullet(Z,\Q),\star)$ denotes the $H$-equivariant Chow groups of $Z$ equipped equipped with the convolution product
$$c\star c^\prime = (p_{13})_*\left(p_{12}^*(c)\cap p_{23}^*(c^\prime)\right),\text{ for } c,c^\prime\in \Chow^\bullet_H(Z,\Q),$$ where $p_{ij}$ denote the projection maps of the triple product $\rnil\times_{\nil}\rnil\times_{\nil}\rnil.$ Using the explicit description of $(\Chow_H^\bullet(Z,\Q),\star)$ (see for example \cite{Chowringsteinberg}) this yields
\begin{corollary} There are equivalences of categories 
\begin{align*}
\MTSpring{G}{\nil}&\cong \Derb(\rmodulesZ\Q[W]\#\Sym^\bullet(X(T)))\text{ and}\\
\MTSpring{G\times \Gm}{\nil}&\cong \Derb(\rmodulesZ\overline{\mathbb{H}}).
\end{align*}
Here $X(T)$ denotes the character group of a maximal torus $T\subset G$, $W$ the Weyl group of $G$ and $\Q[W]\#\Sym^\bullet(X(T))$ the semidirect product of the group algebra of $W$ with the symmetric algebra of $X(T).$ Furthermore $\overline{\mathbb{H}}$ denotes the \emph{graded affine Hecke algebra} associated to $G$ as defined by Lusztig \cite{MR991016}.
\end{corollary}
\subsection{Relation To Other Work} %(1) The first half of this paper strongly borrows from \cite{DLP}, we show that their methods translate to motives.\\
%\noindent (2) 
The second half of this paper is a motivic version of the \emph{derived Springer correspendence} as constructed by Rider \cite{LR1} in the context of equivariant mixed $\ell$-adic sheaves. The motivic setup, as very recently introduced by Soergel--Wendt--Virk \cite{SoeWeVir}, has certain advantages. First, there are no non-trivial extensions between Tate objects, corresponding to the vanishing of rational higher $K$-theory of finite fields and their algebraic closures. Using this, technical difficulties with $\ell$-adic sheaves that necessitate to state the derived Springer correspondence in terms of either dg-categories or the homotopy category of pure complex disappear. Second, Springer motives are defined rationally. Hence all statements are independent of $\ell.$
\subsection{Future Work} (1) In upcoming work with Shane Kelly, generalizing \cite{EKe}, we define a formalism of equivariant motives with coefficients in a finite field. This will allow us to prove a modular motivic derived Springer correspondence analogously. \\\noindent (2) It would be interesting to also consider a generalized motivic derived Springer correspondence along the lines of Rider--Russell \cite{LR2} and \cite{LR3}.\\
(3) Also, it would be interesting to obtain a $K$-theoretic version of the result using equivariant $K$-motives in the sense of Hoyois, see \cite{Hoy1}\cite{Hoy2}. Then it would be possible to handle the affine Hecke algebra. Similar $K$-motivic statements have considered by the author in the case of flag varities \cite{EKo}.
%\subsection{Conventions} 
%By a $\C$-algebra $A$ we always mean a (not necessarily commutative) $\C$-algebra with unit. By $A\modules$ we denote the category of \emph{finitely generated} $A$-modules. If $A=\bigoplus_{n\in\Z}A_n$ is additionally $\Z$-graded, we denote by $A\modules^\Z$ the category of graded $A$-modules and by $A\modules^{\Z,ev}$ the category of evenly graded modules, i.e. those modules which are concentrated in even degrees. For a graded module $M$ its $n$-th shift by $M\langle n\rangle$, where
%$(M\langle n\rangle)^i=M^{i+n}.$

%For an abelian category $\mathcal A$, we denote by $\Der(\mathcal A)$ and $\Derb(\mathcal A)$ its (bounded) derived category and by $\Proje\mathcal A$ the full additive subcategory of projective objects in $\mathcal A$.
%For an additive category $\mathcal A$, we denote by $\Hot(\mathcal A)$ and $\Hotb(\mathcal A)$ its (bounded) homotopy category of chain complexes.
\subsection{Acknowledgements} We thank Rapha\"el Rouquier and Wolfgang Soergel for helpful discussions. We thank George Lusztig for answering a question about \cite{DLP}. Also, we would like to warmly thank the referee for a constructive report.
%Todo:
%\begin{itemize}
%	\item Kapiteleinleitungen
%	\begin{itemize}
%		\item Darstellungstheorie
%		\item Geometrie
%		\item Setup (vorallem bei den Garben und Tilting beschreiben dass das kopiert ist)
%		\item Resultate.
%	\end{itemize}
%	\item ZWEITKORREKTOR OO
%	\item Wie zitiere ich meine eigene ArXiv Arbeit
%\end{itemize}

%% file: Inhalt/MotiveofSpr.tex
%!TEX root = ../main.tex
In this section we show how to translate the results of \cite{DLP} and \cite[Chapter 11]{JanNil} to motives and prove Theorem \ref{thm:pureTate}.
In the following, $k$ denotes an algebraically closed field and $\Lambda$ a commutative ring of coefficients, such that either
\begin{enumerate}
\item resolution of singularities holds over $\Spec(k)$ or
\item the exponential characteristic of $k$ is invertible in $\Lambda.$
\end{enumerate}
For all standard properties of motives we refer to \cite[Sections 14, 16]{MVW} and \cite[Section 5.3]{SKA}. While \cite{MVW} assumes resolution of singularities for many statements about motives of singular schemes, \cite{SKA} shows that requiring that the exponential characteristic of $k$ is invertible in the coefficient ring $\Lambda$ suffices. 
For a variety $X\in\Var_k$ over $k$ we denote its motive by $M(X)\in\DM(\Spec(k), \Lambda)$ and its motive with compact support by $M^c(X)\in\DM(\Spec(k), \Lambda).$ Denote the Tate motive by $\Lambda(1)\in\DM(\Spec(k), \Lambda).$ By definition, the Tate motive is a shift of the reduced motive of $\mathbb{G}_m$, namely the cone of the natural morphism $M(\mathbb{G}_m)[-1]\rightarrow M(\Spec(k))[-1].$
\subsection{Tate motives} We state some general results on pure Tate motives.
\begin{definition} A motive $M\in\DM(\Spec(k), \Lambda)$ is called \emph{pure Tate} if it is isomorphic to a finite direct sum of Tate motives of the form $\Lambda(n)[2n].$
\end{definition}
%\subsection{Properties (MT) and (PT)} Let $X\in \Var_k$ be some variety.
%We introduce the following properties of $X.$
%\begin{itemize}
%\item[(MT)] $M(X)\in \DM(k)$ is \emph{mixed Tate}, that is, $M(X)$ is contained in the triangulated category generated by the Tate motives $\Z(n)\in\DM(k).$
%\item[(PT)] $M(X)\in \DM(k)$ is \emph{pure Tate}, that is, $M(X)$ is a finite direct sum of the Tate motives $\Z(n)[2n]\in\DM(k).$
%\end{itemize}
%Using the properties of $\DM$ (see \cite[Properties 14.5]{MVW}) we immediately get
%\begin{lemma}Let $X\in \Var_k$ and let $E\rightarrow X$ be vector bundle.
%\begin{enumerate}
%\item  Let $X_1,\dots, X_s$ an $\alpha$-partition $X$ such that all $X_i$ have property $(MT).$ Then $X$ has property (MT).
%\item If $X$ admits and affine paving, then $X$ has property (MT).
%\item If $X$ has property (MT), resp. (PT), then so does $E$.
%\item If $X$ has property (MT), resp. (PT), then so does $\Proj(E)$, the projectivization of $E.$
%\end{enumerate}
%\end{lemma}
%\begin{lemma} Let $X\in\Var_k$ be smooth projective. If $X$ has property (MT) then $X$ has property (PT).
%\end{lemma}
%\begin{proof}
%TODO
%\end{proof}
%\begin{lemma} Let $X\in\Var_k$ be smooth projective with an action of a torus $T$ such that $X^T$ hast property $(MT)$ then so does $X.$
%\end{lemma}
%\begin{proof}
%TODO
%\end{proof}
%\begin{lemma} Let $X\in\Var_k$ be smooth quasi-projective variety and $Z\subset X$ a closed subvariety. If $X$ and $Z$ have property (MT), resp. (PT), so does $Y=\Bl_Z(X)$, the Blow-up of $X$ along $Z.$
%\end{lemma}
%\begin{proof}
%TODO
%\end{proof}
%
\begin{definition} Let $X\in \Var_k.$ An \emph{$\alpha$-partition} of $X$ is a finite family of subvarieties  $X_1,\dots, X_s$ of $X$ such that $\bigcup_{i=1}^r X_i$ is closed in $X$ for all $1\leq r\leq s.$ If each $X_i$ is furthermore isomorphic to some affine space $\A^n$, we call the \emph{$\alpha$-partition} an \emph{affine paving} of $X.$
\end{definition}

\begin{lemma} Let $X,Y\in \Var_k$ and $E\rightarrow X$ a vector bundle of rank $r.$\label{lem:puretategeneral}
\begin{enumerate}
\item There is an isomorphism $\Lambda(p)[2p]\otimes\Lambda(q)[2q]\cong \Lambda(p+q)[2(p+q)].$
\item If $X$ is proper, then $M^c(X)=M(X).$
\item\label{lem:puretategeneral:kuen} There is an isomorphism $M(X\times Y)=M(X)\otimes M(Y)$ and if $M(X)$ and $M(Y)$ are pure Tate, then so is $M(X\times Y).$
\item There is an isomorphism $M^c(X\times Y)=M^c(X)\otimes M^c(Y)$ and if $M^c(X)$ and $M^c(Y)$ are pure Tate, then so is $M^c(X\times Y).$
\item There is an isomorphism $M^c(E)\cong M^c(X)(r)[2r]$ and $M^c(X)$ is pure Tate if and only if  $M^c(E)$ is pure Tate.
\item There is an isomorphism $$M^c(\Proj(E))\cong \bigoplus_{p=0}^{r-1}M^c(X)(p)[2p]$$ and $M^c(X)$ is pure Tate if and only if  $M^c(\Proj(E))$ is pure Tate.
\item Let $X$ be smooth and $Z\subset X$ be a closed smooth subvariety of codimension $c$. Denote the the blow-up of $X$ along $Z$ by $\Bl_Z(X).$ Then 
 $$M(\Bl_Z(X))\cong M(X)\oplus\bigoplus_{p=1}^{c-1}M(Z)(p)[2p]$$
 and $M(X)$ and $M(Z)$ are pure Tate if and only if $M(\Bl_Z(X))$ is pure Tate.
\item If $Z\subset X$ is a closed subvariety and $U=X \backslash Z.$ Then there is a distinguished triangle, called localisation triangle,
\begin{center}
\disttriangle{M^c(Z)}{M^c(X)}{M^c(U)}
\end{center} 
and $M^c(X)$ is pure Tate if and only if $M^c(Z)$ and $M^c(U)$ are pure Tate.
\item If $X_1,\dots, X_s$ is an $\alpha$-partition of $X$, then  $M^c(X)$ is pure Tate if and only if $M^c(X_i)$ is pure Tate for all $i.$
\item If $X$ has an affine paving, then $M^c(X)$ is pure Tate.
\end{enumerate}
\end{lemma}
\begin{proof}
(1)-(7) Can be found in \cite[Sections 14-16]{MVW}.
\\ \noindent (8) To show that $M^c(X)$ is pure Tate if $M^c(Z)$ and $M^c(U)$ are, we claim that 
the boundary map $M^c(U)\rightarrow M^c(Z)[1]$ in the localisation triangle vanishes. Let $\Lambda(q)[2q]$ and $\Lambda(p)[2p]$ be direct summands of $M^c(U)$ and $M^c(Z)[1]$, respectively. Then
\begin{align*}
\Hom{\DM(\Spec(k),\Lambda)}{\Lambda(q)[2q]}{\Lambda(p)[2p+1])}
%&=\Hom{\DM(\Spec(k),\Lambda)}(\Lambda(q-p)[2(q-p)-1],\Lambda)\\
&=\Chow^{p-q}(\Spec(k),-1)\otimes \Lambda=0
\end{align*}
where the right hand side denotes a higher Chow group which vanishes since in general $\Chow^\bullet(-,-1)=0.$
Hence $M^c(X)=M^c(Z)\oplus M^c(U)$ and the statement follows.
\\ \noindent (9) Follows by (8) and induction.
\\ \noindent (10) Follows from $M^c(\A^n)=\Lambda(n)[2n]$ and (9).
\end{proof}
As demonstrated in \cite{BBB}, there is a Bia\l{}ynicki-Birula decomposition of motives of varieties with $\Gm$-actions. This can sometimes be used to show that a smooth projective variety is pure Tate.
\begin{lemma}\label{lem:bbstrat}
Let $X\in \Var_k$ be a smooth projective variety equipped with an action of $\Gm.$ Then $M(X^\Gm)=M^c(X^\Gm)$ is pure Tate if and only if $M(X)=M^c(X)$ is.
\end{lemma}
\begin{proof} The Bia\l{}ynicki-Birula theorem gives an $\alpha$-partition of $X$ into vector bundles on the connected components of $X^\Gm.$ The statement follows using the previous lemma.
\end{proof}
\subsection{Prehomogeneous Vector Spaces and Pure Tateness}\label{subsec:prehom}
We show how the methods of \cite[Section 2]{DLP} allow to prove that a certain variety associated to a prehomogeneous vector space is pure Tate. We recall some of their notation. 

Let $M$ be a connected linear algebraic group and $V$ a prehomogeneous $M$-module, that is, $V$ contains a dense $M$-orbit $V^0.$ Fix a $v\in V^0$ and denote by $M_v$ its stabilizer in $M.$ Let $H$ be a closed Borel-subgroup in $M$ and $U$ an $H$-stable linear subspace of $V.$ Let 
\begin{align*}
M_U&=\{g\in M\,|\, g^{-1}v\in U\} \text{ and}\\
X_U&=M_U/H.
\end{align*}
We are interested in the motive of the varieties $X_U.$ By \cite[Lemma 2.2(i)]{DLP} and since $H$ is a Borel subgroup, $X_U$ is a smooth projective variety.

Let $\Gamma$ be the set of $H$-stable subspaces of $V.$ For $U\in\Gamma$, let $P_U$ be the stabilizer of $U$ in $M$ and denote
$\delta(U)=\dim(M/P_U)-\dim(V/U)$ and $\gamma(U)=\dim(M/H)-\dim(V/U).$
In \cite[Section 2.7]{DLP}, $\Gamma$ is equipped with the structure of a directed graph, whose edges $(U, U^\prime)$ have the property
\begin{enumerate}
\item $U\subset U^\prime$,
\item $\dim U^\prime/U=1$ and 
\item there exists a parabolic subgroup $P\supset H$ of semisimple rank 1 and $U^{\prime\prime}\in \Gamma$ such that $U^{\prime\prime}\subset U$, $P\subset P_{U^{\prime\prime}}$, $P\not\subset P_U$, $P\subset P_{U^\prime}$ and $\dim U/U^{\prime\prime}=1.$
\end{enumerate}
 Then, the $M$-module $V$ is called \emph{good} if for any $U\in \Gamma$ either $U\subset U^\prime$ for some $U^\prime\in \Gamma$ with $\delta(U^\prime)<0$, or $U$ lies in the same component of $\Gamma$ as some $U^\prime\in\Gamma$ with $\delta(U^\prime)\leq 0.$
 
 In \cite[Proposition 2.12]{DLP} it is then shown that under the condition that $V$ is good the Borel--Moore-homology of $X_U$ is torsion free, concentrated in even degrees and generated by algebraic cycles. We copy their arguments and show that $M(X_U)$ is pure Tate.
 
 \begin{proposition}\label{prop:goodTate} Assume that the $M$-module $V$ is good. Then $M(X_U)$ is pure Tate for all $U\in \Gamma.$
 \end{proposition}
 \begin{proof}
To translate the inductive argument used in \cite{DLP} to the world of motives, we will use the slice filtration $\nu$ for effective motives as studied in \cite{huka}. 

To any $X_U$, this associates a family of objects $\nu_{<n}M(X_U)$ for $n\geq 0$ and a family of compatible morphisms $\nu_{<n}M(X_U)\rightarrow\nu_{<m}M(X_U)$ for $n\geq m.$ By definition $\nu_{<0}M(X_U)=0$ and since $X_U$ is a smooth projective variety by \cite[Propositions 1.7, 1.8]{huka} we have $\nu_{<n} M(X_U)=M(X_U)$ for $n>\dim X_U$. 

Hence it suffices to show that $\nu_{<n}M(X_U)$ is pure Tate for all $n.$ For each $U\in \Gamma$ and $n\geq 0$ consider the following statement
\begin{enumerate}
\item[$(PT_n)$] $\nu_{<n}M(X_U)$ is pure Tate. % if $i\leq \gamma(U)-n.$
\end{enumerate}
We prove this by induction on $n.$ If $n=0$, then $\nu_{<n}M(X_U)=0.$ 
Now let $n\geq 1$ and assume that $(PT_{n-1})$ holds for all $U\in\Gamma.$ 

If $U\in\Gamma$ is contained in $U^\prime\in\Gamma$ with $\delta(U^\prime)<0$, then $X_U$ is empty and hence $M(X_U)=0$ pure Tate. If $\delta(U)\leq 0$, then $X_U$ is a finite disjoint union of flag varieties isomorphic to $P_U/H$ by \cite[Paragraph 2.9 (a)]{DLP}. The Bruhat decomposition provides an affine paving of $P_U/H.$ Hence  $M(X_U)$ is pure Tate by Lemma \ref{lem:puretategeneral}.

Since $V$ is good, for each connected component of $\Gamma$ there is hence some $U$ for which $(PT_n)$ holds. So the statement of the proposition reduces to the following lemma.
  \end{proof}
\begin{lemma}
Assume that $n\geq 1$ and $(PT_{n-1})$ holds for all $U_1\in\Gamma.$ Let $U\subset U^\prime$ be an edge in $\Gamma.$ Then  $(PT_n)$ holds for $U$ if and only if $(PT_n)$ holds for $U^\prime.$
\end{lemma}
\begin{proof}
Let $U^{\prime\prime}$ and $P$ as in property (3) of edges of $\Gamma.$
Let 
$$Z =\{(gH, g^\prime H) \in  M/H\times M/H \,|\, g^{-1}v\in U \text{ and } g^{-1}g^\prime\in P\}.$$
Then by \cite[Lemma 2.11]{DLP} $Z\rightarrow X_U$ is the projectivization of a vector bundle of rank two $E\rightarrow X_U$ and furthermore $Z\cong \Bl
_{X_{U^{\prime\prime}}}(X_{U^\prime})$, where $X_{U^{\prime\prime}}\subset X_{U^\prime}$ is a closed subvariety of codimension two. Hence by the projective bundle and blow-up formula we have
$$
M(X_U)\oplus M(X_U)(1)[2]=M(Z)=M(X_{U^\prime})\oplus M(X_{U^{\prime\prime}})(1)[2].
$$
Applying $\nu_{<n}$ yields that
$$ \nu_{<n}\left(M(X_U)\oplus M(X_U)(1)[2]\right)=\nu_{<n}(M(X_U))\oplus \nu_{<n-1}(M(X_U))(1)[2]$$
equals
$$ \nu_{<n}\left(M(X_{U^\prime})\oplus M(X_{U^{\prime\prime}})(1)[2]\right)=\nu_{<n}(M(X_{U^\prime}))\oplus \nu_{<n-1}(M(X_{U^{\prime\prime}}))(1)[2]$$
where we use that $\nu_{<n}(-(1))=\nu_{<n-1}(-)(1)$, see \cite[Corollary 1.4(v)]{huka}.
Now $\nu_{<n-1}(M(X_{U^{\prime\prime}}))$ and $\nu_{<n-1}(M(X_U))$ are pure Tate by induction, and the Statement follows.
\end{proof}
\subsection{Springer Fiber is Pure Tate} We prove Theorem \ref{thm:pureTate} from the introduction. Let $G$ be reductive algebraic group over $k$, $N\in \nil\subset \g$ be an element of the nilpotent cone in the Lie algebra of $G$ and $\B{N}=\mu\inv(N)\subset\rnil$ the Springer fiber in the Springer resolution $\mu: \rnil\rightarrow\nil.$  Let $$\B{}=\{\borel\subset \g \,|\, \borel \text{ is a Borel subalgebra}\}$$ denote the flag variety. Then we can identify
$$\B{N}=\{\borel\subset \g \,|\, N\in \borel\}\subset\B{}.$$
The goal is to prove that $M(\B{N})\in \DM(\Spec(k),\Lambda)$ is pure Tate. We note that the Springer fiber is proper and hence $M^c(\B{N})=M(\B{N}).$

As the Springer fiber only depends on the isogeny class of the semisimple part of $G,$ we may assume that $G$ is of adjoint type and hence a direct product of its simple constituents (see \cite[Section 2.7]{JanNil}). 
Furthermore, a Springer fiber of a direct product of groups decomposes into a direct product as well, and by Lemma \ref{lem:puretategeneral}(\ref{lem:puretategeneral:kuen}) it suffices to consider each individual factor. 

So we can assume that $G$ is a simple algebraic group. If $G$ is a classical group, that is, of type $A,B,C$ or $D,$ then $\B{N}$ admits a paving by affine spaces  by \cite{DLP} if $\chara(k) =0$ and more general by \cite[Theorem 11.22]{JanNil} if $p=\chara(k)$ is good for $G.$ Hence $M(\B{N})$ is pure Tate in this case by Lemma \ref{lem:puretategeneral}.

We can hence assume that $G$ is a simple group of exceptional type $E, F$ or $G$ and assume that $p>3(h+1)$, where $h$ denotes the Coxeter-number of $G.$ 
%We want to further simplify to the case that $N\in \nil$ is a distinguished nilpotent element. For this 
We proceed as in \cite[Section 3.4]{DLP}.

There exists a special cocharacter $\tau: \Gm\rightarrow G$ associated to $N$, see  \cite[Section 5.2]{JanNil} for a definition and existence and uniqueness result, alternatively use the Morozov--Jacobson theorem, which holds since $p>3(h-1)$. This cocharacter induces a decomposition of $\g$ into even weight spaces
 $$\g=\bigoplus_{n\in\Z} \g(2n, \tau),$$
 such that $N\in \g(2,\tau).$ Let $G_0\subset P\subset G$ be the Levi and parabolic subgroup with Lie algebra $\g(0,\tau)$ and $\bigoplus_{n\geq 0}\g(2n, \tau)$, respectively, and let $S=\tau(\Gm).$
Now $\B{N}$ admits an $\alpha$-filtration by intersecting it with the $P$-orbits $\mathcal O$ on $\B{}.$ Each of those intersections $\B{N,\mathcal O}=\B{N}\cap \mathcal O$ is smooth projective. 

So $\M(\B{N})$ is pure Tate if and only if $M(\B{N,\mathcal O})$ is pure Tate for each $\mathcal{O}$ by Lemma \ref{lem:puretategeneral}.
Furthermore each  $M(\B{N,\mathcal O})$ is pure Tate if and only if $M(\B{N,\mathcal O}^S)$ is pure Tate by Lemma \ref{lem:bbstrat}.

By a similar argument to \cite[Section 3.6]{DLP} and using  Lemmata \ref{lem:puretategeneral}, \ref{lem:bbstrat} again, we can reduce to the case that $N$ is in fact a distinguished nilpotent element, so not already contained in the Lie algebra of any proper Levi subgroup of $G.$

Let $B_0\subset G_0$ be a Borel subgroup. Denote its Lie algebra by $\borel_0.$ Then there is a unique $\borel_\mathcal{O}\in \mathcal{O}^S$ with $\borel_\mathcal{O}\cap \g(0,\tau)=\borel_0.$ Let $U_\mathcal{O}=\borel_\mathcal{O}\cap \g(2,\tau).$ This is a $B_0$-stable linear subspace of the prehomogeneous $G_0$-module $\g(2,\tau).$ We can hence consider the variety $X_{U_\mathcal{O}}$ as defined in Section \ref{subsec:prehom}. In fact the map $X_{U_\mathcal{O}}\rightarrow \B{N,\mathcal O}^S, \,gB_0\mapsto g\borel_\mathcal{O}$ is an isomorphism.

Now \cite{DLP} show by an involved case by case computation that the $G_0$-modules $\g(2,\tau)$ arising in this way from a distinguished nilpotent element for an exceptional group are \emph{good.} This computation, as the whole paper, is carried out for $k=\C$ but also works as long as the Morozov--Jacobson theorem holds, so in particular if $p>3(h-1).$ We thank George Lusztig for answering a question about that. We can hence use Proposition \ref{prop:goodTate} to see that $M(X_{U_\mathcal{O}})=M(\B{N,\mathcal O}^S)$ is pure Tate.

This concludes the proof of Theorem \ref{thm:pureTate}.

%Denote by $G^0_N\subset G$ the connected component of the identity in the stabilizer of $N$ in $G.$ Let $S^\prime\subset G^0_X$ be a maximal torus. Then the centralizer $L^\prime=\operatorname{Z}_G(S^\prime)$ of $S^\prime$ in $G$ is a Levi subgroup of $G.$ Then $N$ is called \emph{distinguished} if $G=L^\prime.$ Furthermore $N\in\mathfrak{l}^\prime=\operatorname{Lie}(L^\prime)$ is distinguished. 
%Denote the flag variety of $L^\prime$ by $\B{}(L^\prime).$

%% file: Inhalt/Form.tex
%!TEX root = ../main.tex
In this section we prove Theorem \ref{thm:motdersprin}. We assume that $k=\overline{\mathbb{F}}_p$ and $\Lambda=\Q$. We denote a variety $X\in\Var_k$ with an action of a linear group $H$ by $(H\looparrowright X)\in\Var_k.$ Morphism between varieties with action are given by pairs 
$$(\phi,f): (H_1\looparrowright X_1)\rightarrow(H_2\looparrowright X_2)$$
of a morphism of linear groups and a morphism of varieties compatible with the actions. If $\phi=\id$ is the identity morphism, we will often drop it from the notation.
Now \cite[Chapter I]{SoeWeVir} associates to the the datum $(H\looparrowright X)$
 the $\Q$-linear tensor triangulated category\footnote{We work with $\Deri(-)=\DA(-,\Q),$ the homotopical stable algebraic derivator of \'etale motives with rational coefficients over the category of varieties $\Var_k$ over $k.$ We note that $\DA(\Spec(k),\Q)\cong \DM(\Spec(k),\Q)$ which is the category of motives we considered in the first part of the paper.} of \emph{$H$-equivariant $\Deri$-motives on $X$} denoted by $\DerG{H}{X}$. We denote the tensor unit by $\un=\un_X$. The system of categories $\DerG{-}{-}$ comes equipped with a six-functor-formalism and induction/restriction functors, very similar to the equivariant derived category of Bernstein--Lunts \cite{BeLu}.
In the case that $X$ has finitely many $H$-orbits, \cite[Chapter II]{SoeWeVir} defines the category \emph{$H$-equivariant orbitwise mixed Tate motives} $\MTDer{H}{X}\subset \DerG{H}{X}$ which are analogous to constructible equivariant sheaves. 
From now on we consider the case $H=G$ or $H=G\times \Gm$ and $X=\nil.$
\subsection{Orbitwise Pure Tateness of the Springer motive} In the introduction we cheated a bit. A priori, it is not clear that the Springer motive $\mu_!(\un_{\rnil})=\mu_*(\un_{\rnil})\in \DerG{H}{\nil}$ already lives in the subcategory $\MTDer{H}{\nil}.$ In this section we show how the pure Tateness of the Springer fiber implies this and that $\mu_!(\un_{\rnil})$ is additionally pointwise pure.
\begin{theorem} Let $\mathbb{O}$ be an $H$-orbit on $\nil$. Let $N$ be a point in $\mathbb{O}$ and denote by $H_N\subset H$ its stabilizer.
Denote the corresponding morphisms of varieties with group action by
\begin{align*}
\mu=(\id,\mu):(H\looparrowright \rnil)&\rightarrow (H\looparrowright\nil),\\
j=(\id,j):(H\hookrightarrow\mathbb{O})&\hookrightarrow (H\looparrowright\nil)\text{ and}\\
(\iota,i):(H_N\looparrowright \{N\})&\hookrightarrow (H\looparrowright\mathbb{O}).
\end{align*}
Then for $?\in\{*,!\}$ we have a chain of functors
\begin{center}
\begin{tikzcd}
\DerG{H}{\rnil}\arrow{r}{\mu_!=\mu_*}&\DerG{H}{\nil}\arrow{r}{j^?}&\DerG{H}{\mathbb{O}}\arrow{r}{(\iota,i)^*}&\DerG{H_N}{\{N\}}\arrow{r}{\For}&\DerG{}{\Spec(k)}
\end{tikzcd}
\end{center}
where $\For$ is the functor forgetting the $H_N$ action.
We claim that 
$$\For(\iota,i)^*j^?\mu_!(\un_{\rnil})\in\DerG{}{\Spec(k)}$$
is pure Tate, that is, a finite direct sum of Tate motives $\un(n)[2n].$
\end{theorem}
\begin{proof} Since the forgetful functor $\For$ commutes with all six functors, it suffices to show the corresponding statement where we forget about all group actions. Since $\mu$ is proper and hence $\mu_*=\mu_!$, it suffices to show the statement for $?=*$ by duality. Now we apply base change for the cartesian diagram
\begin{center}
\begin{tikzcd} \B{N}=\mu^{-1}(N) \arrow{r}{k}  \arrow{d}{\fin_{\B{N}}}   &  \mu^{-1}(\mathbb{O}) \arrow{r}{l} \arrow{d}{}& \rnil \arrow{d}{\mu} \\ 

\Spec(k)=\{N\}  \arrow{r}{i} &  \mathbb{O} \arrow{r}{j} & \nil
	\end{tikzcd}
%\cartesiandiagramwithmaps{\B{N}}{k}{\rnil}{\mu_N}{\mu}{\{N\}}{l}{\nil}
\end{center}
and see that 
  $$\For(\iota,i)^*j^*\mu_!(\un_{\rnil})\cong\fin_{\B{N},!}(k^*l^*\un_{\rnil})=\fin_{\B{N},!}(\un_{\B{N}})\in\DerG{}{\Spec(k)}.$$
Since we are working with rational coefficients there is a natural equivalence $$\DerG{}{\Spec(k)}=\DA(\Spec(k),\Q)\cong \DM(\Spec(k),\Q),$$ see \cite{Ayoubetale}, and the Verdier dual of $\fin_{\B{N},!}(\un_{\B{N}})\in\DerG{}{\Spec(k)}$ corresponds to  $M^c(\B{N})\in\DM(\Spec(k),\Q)$, the motive of the Springer fiber. But $M^c(\B{N})$ is pure Tate by Theorem \ref{thm:pureTate}.
\end{proof}
 By definition, $\MTDer{H}{\nil}$ consists of motives $M$ such that $$\For(\iota,i)^*j^?M\in\left\langle \un(n) \,|\,n\in\Z\right\rangle_\Delta\subset\DerG{}{\Spec(k)}$$ is mixed Tate, i.e. contained in the triangulated category generated by motives $\un(n)\in \DerG{}{\Spec(k)},$ for each orbit $j:\mathbb{O}\hookrightarrow \nil, ?\in\{*,!\}$ and $(\iota,i)^*$ defined as above. A motive $M$ is called \emph{pointwise $?$-pure} if additionally  $\For(\iota,i)^*j^?M$ is pure Tate  (so a finite direct sum of motives $\un(n)[2n]$) for $?\in\{*,!\}$ and \emph{pointwise pure} if it is pointwise $*$-pure and pointwise $!$-pure.
\begin{corollary} We have $\mu_*(\un_{\rnil})\in\MTDer{H}{\nil}$ and $\mu_*(\un_{\rnil})$ is pointwise pure.
\end{corollary}
\subsection{Tilting and Formality for Springer motives}
Denote by
$$\MTSpring{H}{\nil}_{add}=\left\langle \mu_*(\un_{\rnil})(n) \,|\, n\in \Z\right\rangle_{\oplus, \inplus}\subset \MTSpring{H}{\nil}$$
the idempotent closed additive subcategory of $\MTSpring{H}{\nil}$ generated by the Springer motive and its Tate twists. 
\begin{lemma}\label{lem:springtilt}
The category $\MTSpring{H}{\nil}_{add}$ is a \emph{tilting subcategory} of $\DerG{H}{\nil},$ meaning that it is
\begin{enumerate}
\item  idempotent closed, additive and
\item $\Hom{\DerG{H}{\nil}}{M}{N[i]}=0$ for all $M,N\in  \MTSpring{H}{\nil}_{add}$ and $i\neq 0.$
\end{enumerate}
  \end{lemma} 
  \begin{proof} (1) Holds by construction.\\\noindent (2) This is implied by the pointwise purity of the objects in $\MTDer{H}{\nil}_{add}$ using an argument very similar to \cite[Corollary 6.3]{SoeWe}.  More, generally we show that for all $M,N\in\MTSpring{H}{\nil}$ with $M$ pointwise $*$-pure and $N$ pointwise $!$-pure $$\Hom{\DerG{H}{\nil}}{M}{N[i]}=0$$ for $i\neq 0.$ For this, denote by $j: U\hookrightarrow \nil$ and $i:Z=\nil\backslash U\hookrightarrow \nil$ the inclusion of the open orbit and its closed complement. Then the localisation sequence
 \begin{center}
\disttriangle{j_!j^!M}{M}{i_*i^*M}
\end{center} 
yields an exact sequence
\begin{center}
\begin{tikzcd}[column sep=small]
 \Hom{\DerG{H}{Z}}{i^*M}{i^!N[i]}\arrow[r] & \Hom{\DerG{H}{\nil}}{M}{N[i]} \arrow[r] & \Hom{\DerG{H}{U}}{j^*M}{j^!N[i]}.
\end{tikzcd}
\end{center}
Now $i^*M$ is pointwise $*$-pure and $i^!N$ is pointwise $!$-pure. Hence, the the left hand side vanishes for $i\neq 0$ by induction (here we use that there are only finitely many orbits). Also, the right hand side vanishes for $i\neq 0$ using the purity assumption, the induction equivalence and that $\Hom{\DerG{H_N}{\nil}}{\Q}{\Q(p)[q]}=0$ unless $p=2q,$ where $H$ denotes the stabiliser of a point $N\in U.$
  \end{proof}
  \begin{remark} The arguments in the proof of Lemma \ref{lem:springtilt} (2) argument do not translate to the non-equivariant case immediately. Hence, it is not clear if a non-equivariant analogue of the motivic derived Springer correspondence is true. One would need to show that the higher Chow groups of nilpotent orbits vanish, which is not clear to the author.
  \end{remark}
  We can now apply the tilting formalism and show
  \begin{theorem}
  There is an equivalence of categories, called \emph{tilting}, between
  $$\tilt: \Hotb(\MTSpring{H}{\nil}_{add})\stackrel{\sim}{\rightarrow}\MTSpring{H}{\nil}$$
  the homotopy category of bounded chain complexes $\MTSpring{H}{\nil}_{add}$ and the category of Springer motives $\MTSpring{H}{\nil}.$
  \end{theorem}
  \begin{proof}
  Since by Lemma \ref{lem:springtilt} $\MTSpring{H}{\nil}_{add}$ is a tilting subcategory of $\DerG{H}{\nil}$ by \cite[Theorem B.3.1]{SoeWeVir}, which establishes a tilting formalism for stable derivators, there is a fully faithful functor
  $$\tilt: \Hotb(\MTSpring{H}{\nil}_{add})\rightarrow\DerG{H}{\nil}.$$
  The essential image of $\tilt$ is the triangulated subcategory of $\DerG{H}{\nil}$ generated by $\MTSpring{H}{\nil}_{add}$ which is $\MTSpring{H}{\nil}$ by definition.
 \end{proof}
We name the $\Z$-graded ``Ext"-algebra\footnote{We put the quotation marks here because of the Tate twists which is not part of the definition of the true Ext-algebra. In fact, the true Ext-algebra is concentrated in degree 0, by Lemma  \ref{lem:springtilt}.} of the Springer motive $\mu_*\left(\un_{\rnil}\right)$
 $$E:=\bigoplus_{n\in\Z} E^n\text{, where }E^n=\Hom{\DerG{H}{\nil}}{\mu_*\left(\un_{\rnil}\right)}{\mu_*\left(\un_{\rnil}\right)(n)[2n]}.$$
 We denote the shift of grading functor for graded $E$-modules by $\langle-\rangle.$
 We can rephrase the last theorem as
 \begin{corollary}
 There is an equivalence of categories
 $$\Derb(\rmodulesZ E)\rightarrow\MTSpring{H}{\nil}$$
 between the bounded derived category of finitely generated graded right modules over $E$ and the category of Springer motives.
 \end{corollary}
 \begin{proof} 
 The category $\MTSpring{H}{\nil}_{add}$ can be identified with the full subcategory
 $$\left\langle E\langle n\rangle \,|\, n\in \Z\right\rangle_{\oplus, \inplus}\subset \rmodulesZ E.$$ 
 of graded right $E$-modules generated by finite direct sums and summand of shifts of $E.$
 By the explicit description of $E$, see the next section, $E$ has finite cohomological dimension, and hence the homotopy category of bounded chain complexes in $\left\langle E\langle n\rangle \,|\, n\in \Z\right\rangle_{\oplus, \inplus}$ is equivalent to the bounded derived category of graded $E$-modules.
 \end{proof}
\subsection{Description of the ``Ext"-algebra $E$} The last step in the proof of Theorem \ref{thm:motdersprin} is to explicitely describe the ``Ext"-algebra $E.$ Recall that $Z=\rnil\times_\nil\rnil$ denotes the Steinberg variety and that $H$ acts on $Z$ by the diagonal action. We repeat some results from \cite[Section 8.6]{crisginz}, who give an explicit description of the sheaf-theoretic version of $E$ in terms of the Borel--Moore homology of the Steinberg-variety. Their results translate to the motivic setting almost word by word, all that is needed is a six-functor formalism.
\begin{lemma}
There is an isomorphism of graded algebras
$$E\cong(\Chow^{\bullet}_H(Z,\Q),\star)$$
between $E$ and the $H$-equivariant Chow groups of $Z$ equipped with the convolution product.
\end{lemma}
\begin{proof}
We just show that $E^n\cong \Chow^{n}_H(Z,\Q)$ as a vector space here. The statements about the algebra structure can be deduced as in \cite[Section 8.6]{crisginz}.
Consider the cartesian diagram of $H$-varieties:
\begin{center}
\cartesiandiagramwithmaps{Z}{p_2}{\rnil}{p_1}{\mu}{\rnil}{\mu}{\nil}
\end{center}
Denote for a variety $X\in\Var_k$ its structure map by $\fin_X:X\rightarrow \Spec(k)=pt.$ Since $\rnil$ is smooth we have $\un_{\rnil}=\fin_{\rnil}^*(\un_{pt})=\fin_{\rnil}^!(\un_{pt})(-d)[-2d]$ where  $d=\dim(\rnil).$ Furthermore note that $\mu$ is proper and hence $\mu_!=\mu_*$. Now consider
\begin{align*}
E^n&=\Hom{\DerG{H}{\nil}}{\mu_!(\un_{\rnil})}{\mu_*(\un_{\rnil})(n)[2n]}\\
%&=\Hom{\DerG{H}{\nil}}{\un_\nil }{\iHom{\nil}(\mu_!\left(\un_{\rnil}\right),\mu_!\left(\un_{\rnil}\right)(n)[2n])}
&=\Hom{\DerG{H}{\rnil}}{\un_{\rnil}}{\mu^!\mu_*(\un_{\rnil})(n)[2n]}\\
&=\Hom{\DerG{H}{\rnil}}{\un_{\rnil}}{p_{1,*}p_2^!(\un_{\rnil})(n)[2n]}\\
&=\Hom{\DerG{H}{\rnil}}{\fin_{\rnil}^*(\un_{pt})}{p_{1,*}p_2^!\fin_{\rnil}^!(\un_{pt})(n-d)[2(n-d)]}\\
&=\Hom{\DerG{H}{pt}}{\un_{pt}}{\fin_{Z,*}\fin_{Z}^!(\un_{pt})(n-d)[2(n-d)]}\\
&=\Chow^{\dim(Z)+(n-d)}_H(Z,\Q)=\Chow^{n}_H(Z,\Q)
\end{align*}
where in the last equality we used \cite[Theorem II.2.9.]{SoeWeVir} (there is no resolution of singularities for $\Spec(k)$ necessary here, using \cite[Section 5.2]{SKA} and that $p$ is invertible in $\Q$) and that $Z$ is of dimension $\dim(Z)=d.$
\end{proof}